\documentclass[10pt]{amsart}
\usepackage[utf8]{inputenc}
\usepackage{amsmath,amsthm,amssymb,amsfonts,mathrsfs,bbm,comment,import,xcolor, hyperref,graphicx}

\newcommand{\der}{\partial}

\newcommand{\g}{\mathfrak{g}}
\newcommand{\A}{\mathcal{A}}
\newcommand{\psih}{\hat{\psi}}
\newcommand{\psit}{\tilde{\psi}}
\newcommand{\One}{\mathbbm{1}}
\newcommand{\QuatSpace}{\mathbb{H}}
\newcommand{\gen}[1]{\ensuremath{\langle #1\rangle}}
\newcommand{\Mat}[1]{\text{Mat}_{#1}(\C)}

\newcommand{\GL}[1]{\text{GL}_{#1}(\C)}

\newcommand{\su}[1]{\mathfrak{su}(#1)}
\newcommand{\hconj}{\dagger}
\newcommand{\Der}{\operatorname{Der}}
\newcommand{\id}[1]{\operatorname{id}_{#1}}
\newcommand{\udots}{\reflectbox{$\ddots$}}

\newcommand{\R}{\mathbb{R}}
\newcommand{\C}{\mathbb{C}}
\newcommand{\Rstar}{\R\setminus\{0\}}

\newcommand{\paraa}[1]{\big(#1\big)}

\theoremstyle{definition}
\newtheorem{definition}{Definition}[section]
\newtheorem{lemma}[definition]{Lemma}
\newtheorem{proposition}[definition]{Proposition}
\newtheorem{theorem}[definition]{Theorem}
\newtheorem{corollary}[definition]{Corollary}
\newtheorem{example}[definition]{Example}

\title{Projective Real Calculi over Matrix Algebras}
\author{Axel Tiger Norkvist}

\begin{document}

\begin{abstract}
    In analogy with the geometric situation, we study real calculi over projective modules and show that they can be realized as projections of free real calculi. Moreover, we consider real calculi over matrix algebras and discuss several aspects of the classification problem for real calculi in this case, leading to the concept of quasi-equivalence of matrix representations. We also use matrix algebras to give concrete examples of real calculi where the module is projective, and show that the existence of a Levi-Civita connection depends on the eigenvectors of specific anti-hermitian matrices in this case. 
\end{abstract}

\maketitle

\section{Introduction}\label{sec:intro}
    \noindent
In the rapidly developing and conceptually rich field of study that is noncommutative geometry it is important to discuss the subject from several points of view in order to understand strenghts and weaknesses of different perspectives. For instance, while spectral triples enable the use of powerful analytic tools to study the properties of noncommutative spaces it can be somewhat difficult to see whether the metric, implicitly given by the Dirac operator, can be realized as a bilinear form over a module. Therefore it can be helpful to take an alternative approach and study the subject from a more direct point of view, where both the metric and the module on which it is realized as a bilinear form are explicitly given, much in the spirit of classical Riemannian geometry. 

In \cite{aw:cgb.sphere, aw:curvature.three.sphere}, the framework of pseudo-Riemannian calculi was introduced in order to discuss the existence (and uniqueness) of a Levi-Civita connection as well as curvature in noncommutative spaces. Although the existence of a Levi-Civita connection cannot be guaranteed in the context of real calculi, it was shown that it is unique whenever it exists. Note that the framework of real calculi is not the only approach to metrics and Levi-Civita connections in noncommutative geometry, and there are several others that are different but share some conceptual similarities (see e.g. \cite{lm:lie.extensions,dvm:connections.central.bimodules,fgr:supersymmetric.noncommutative, bm:starCompatibleConnections,r:leviCivita,mw:quantum.koszul,a:2020cartan.and.LC.connection}).

In \cite{aw:curvature.three.sphere} the Riemannian curvature of the noncommutative 3-sphere was studied, and in \cite{aw:cgb.sphere} a Gauss-Chern-Bonnet theorem was formulated for the noncommutative 4-sphere. As an example of a noncompact noncommutative manifold with nontrivial features, the noncommutative cylinder was featured in \cite{al:noncommutative.cylinder} and a pseudo-Riemannian calculus was constructed in order to explicitly calculate the Levi-Civita connection and the corresponding curvature. In \cite{atn:nc.minimal.embeddings} algebraic aspects of real calculi were discussed with the introduction of real calculus homomorphisms, and these were interpreted geometrically to develop a noncommutative theory of submanifolds. In general, it is interesting to see examples of real calculi where the module used over the algebra to form a real calculus is projective, but in all of the examples given in the above-mentioned articles the modules are even more well-behaved in the sense that they are free. In this paper we shall consider real calculi over projective modules, and in analogy with projective modules being projections of free modules we show that projective real calculi can be realized as "projections" of free real calculi.

Although the study of projective real calculi is an important part of this text, the main focus will lie on real calculi over $\A=\Mat{N}$ and how isomorphisms of real calculi are characterized in this case; in particular, we give a complete classification of real calculi over $\A=\Mat{N}$ in the case where the Lie algebra $\g$ acting on $\A$ as a set of derivations is 1-dimensional and the module is $\C^N$. Real calculi over this algebra are also used to give concrete examples of real calculi where the associated $\A$-module $M=(\C^{N})^n$ is projective (but not free whenever $n\neq N$). From a more general perspective, the study of finite-dimensional algebras is not only done for the purpose of highlighting interesting aspects of noncommutative geometry (for instance, a relatively rich differential geometric structure using $\Omega_D(\Mat{N})$ as differential forms can be developed on matrix algebras despite the fact that derivations are inner, as described in \cite{dvkm:nc.diff.geom.matrix.algebras}), but these algebras have an important role in their own right from a theoretical point of view (see for instance the spectral Standard Model described in \cite{cc:why.nc.standard.model}, where the finite-dimensional algebra $\A_F=\C\oplus\QuatSpace \oplus\Mat{3}$ encodes the symmetries of the standard model). Hence, it is important to study how real calculi over matrix algebras behave, in order to understand how the framework of real calculi fits into the larger body of work in noncommutative geometry as a whole.

The paper is organized as follows. First, in Section \ref{sec:prelims} we recall basic definitions and results regarding real calculi and their morphisms. Then, in Section \ref{sec:free.and.projective} we introduce the concept of projective real calculi, and discuss how these structures can be understood as projections of free real calculi. We describe affine connections on the projective module and how they are derived from a free module before we turn our attention to metrics and projective real calculi. In Section \ref{sec:matrix.algebra} we treat isomorphisms of real calculi over $\A=\Mat{N}$, both when the module $M$ is free and when $M=(\C^N)^n$. Finally, in Sections \ref{subsec:ex.dim.1} and \ref{subsec:ex.dim.n} we show that the existence of a Levi-Civita connection when $\g$ is abelian and $M=(\C^N)^n$ depends on the map $\varphi:\g\rightarrow (\C^N)^n$ and the eigenvectors of specific traceless and anti-hermitian matrices representing hermitian derivations.
\section{Preliminaries}\label{sec:prelims}
    We begin by recalling the basic definitions and results regarding real calculi (see \cite{aw:curvature.three.sphere}) and their morphisms (see \cite{atn:nc.minimal.embeddings}) which make out the framework used throughout this text. The framework takes a derivation-based approach where a module over the algebra representing the (noncommutative) space in question acts as a generalization of the concept of vector fields. However, since the set $\Der(\A)$ generally lacks a module structure when $\A$ is noncommutative the correspondence between vector fields and derivations in the framework is given as an explicit linear map which is not assumed to be an isomorphism. We make the notion of a real calculus precise below.

\begin{definition}[Real calculus]
	Let $\A$ be a unital $^*$-algebra, and let $\g_D$ denote a real Lie algebra together with a faithful representation $D:\g\rightarrow\text{Der}(\A)$ such that $D(\der)$ is a hermitian derivation for all $\der\in\g$. Moreover, let $M$ be a right $\A$-module and let $\varphi:\g\rightarrow M$ be a $\R$-linear map such that $M$ is generated by $\varphi(\g)$. Then
	$C_{\A}=(\A,\g_D,M,\varphi)$ is called a \textit{real calculus} over $\A$.
\end{definition}
\noindent
In the following we will write $\der(a)$ instead of the more cumbersome $D(\der)(a)$ for an element $\der\in\g$ and $a\in\A$ as long as there is no risk of confusion. Moreover, if $\g\subseteq \Der(\A)$ and the representation $D$ is left unspecified this is to be interpreted as $D$ being the identity map.
\begin{definition}
    Let $\A$ be a $^*$-algebra and let $M$ be a right $\A$-module. A \textit{hermitian form} on $M$ is a map $h:M\times M\rightarrow\A$ with the following properties:
    \begin{itemize}
        \item[$h1$.] $h(m_1,m_2+m_3)=h(m_1,m_2)+h(m_1,m_3)$
        \item[$h2$.] $h(m_1,m_2a)=h(m_1,m_2)a$
        \item[$h3$.] $h(m_1,m_2)=h(m_2,m_1)^*$
    \end{itemize}
    for every $m_1,m_2,m_3\in M$ and $a\in\A$. Moreover, if $h(m_1,m_2)=0$ for every $m_2\in M$ implies that $m_1=0$ then $h$ is said to be \textit{nondegenerate}, and in this case we say that $h$ is a \textit{metric} on $M$. The pair $(M,h)$ is called a \textit{(right) hermitian} $\A$-\textit{module}, and if $h$ is a metric on $M$ we say that $(M,h)$ is a \textit{(right) metric} $\A$-\textit{module}.
    Finally, we say that $h$ is \emph{invertible} if the map $\hat{h}:M\rightarrow M^*$, defined by $\hat{h}(m)(n)=h(m,n)$, is a bijection.
\end{definition}

\begin{definition}[Real metric calculus]
    Let $C_{\A}=(\A,\g_D,M,\varphi)$ be a real calculus over $\A$ and let $(M,h)$ be a (right) metric $\A$-module. If
    \begin{equation*}
        h(\varphi(\der_1),\varphi(\der_2))^*=h(\varphi(\der_1),\varphi(\der_2))
    \end{equation*}
    for every $\der_1,\der_2\in\g$ then the pair $(C_{\A},h)$ is called a \textit{real metric calculus}.
\end{definition}
\noindent
Since free modules are the simplest kind of modules, real calculi $C_{\A}=(\A,\g_D,M,\varphi)$ where the module $M$ is free are interesting to study. However, the map $\varphi$ can be degenerate, leading to complications in the analysis of real calculi over a free module. Therefore we shall state additional regularity conditions that $\varphi$ needs to satisfy before the real calculus $C_{\A}$ can be considered free.

\begin{definition}[Free real (metric) calculus]
    Let $C_{\A}=(\A,\g_D,M,\varphi)$ be a real calculus and suppose that $M$ is a free module. If there exists a basis $\{\der_1,...,\der_n\}$ of $\g$ such that $\{\varphi(\der_1),...,\varphi(\der_n)\}$ is a basis of $M$, then $C_{\A}$ is a \textit{free real calculus}. Moreover, if $(C_{\A},h)$ is a real metric calculus and $h$ is an invertible metric, then $(C_{\A},h)$ is a \textit{free real metric calculus}.
\end{definition}

\begin{definition}
    Let $C_{\A}=(\A,\g_D,M,\varphi)$ be a real calculus over $\A$. An \textit{affine connection} on $(M,\g)$ is a map $\nabla:\g\times M\rightarrow M$ satisfying
    \begin{enumerate}
        \item $\nabla_{\der}(m+n)=\nabla_{\der}m+\nabla_{\der}n$
        \item $\nabla_{\lambda\der+\der'}m=\lambda\nabla_{\der}m+\nabla_{\der'}m$
        \item $\nabla_{\der}(ma)=(\nabla_{\der}m)a+m\der(a)$
    \end{enumerate}
    for $m,n\in M$, $\der,\der'\in\g$, $a\in\A$ and $\lambda\in\R$.
\end{definition}

\begin{definition}[Real connection calculus]
    Let $(C_{\A},h)$ be a real metric calculus and let $\nabla$ be an affine connection on $(M,\g)$. Then $(C_{\A},h,\nabla)$ is called a \textit{real connection calculus} if
    \begin{equation*}
        h(\nabla_{\der}\varphi(\der_1),\varphi(\der_2))=h(\nabla_{\der}\varphi(\der_1),\varphi(\der_2))^*
    \end{equation*}
    for every $\der,\der_1,\der_2\in\g$.
\end{definition}

\begin{definition}[Pseudo-Riemannian calculus]
    Let $(C_{\A},h,\nabla)$ be a real connection calculus. We say that $(C_{\A},h,\nabla)$ is \textit{metric} if
    \begin{equation*}
        \der(h(m,n))=h(\nabla_{\der}m,n)+h(m,\nabla_{\der}n)
    \end{equation*}
    for every $\der\in\g$ and $m,n\in M$, and \textit{torsion-free} if
    \begin{equation*}
        \nabla_{\der_1}\varphi(\der_2)-\nabla_{\der_2}\varphi(\der_1)-\varphi([\der_1,\der_2])=0
    \end{equation*}
    for every $\der_1,\der_2\in\g$. A metric and torsion-free real connection calculus is called a \textit{pseudo-Riemannian calculus}.
\end{definition}

A connection fulfilling the conditions for a pseudo-Riemannian calculus is called a \textit{Levi-Civita connection}. In the quite general setup of real metric calculi the existence of a Levi-Civita connection can not be guaranteed. However, it is unique if it exists and if we make the extra assumption that the real metric calculus is free existence can be guaranteed as well.

\begin{theorem}[Levi-Civita connection]\label{thm:existence.and.uniqueness.of.LC}
Let $(C_{\A},h)$ be a real metric calculus. Then there exists at most one affine connection $\nabla$ such that $(C_{\A},h,\nabla)$ is a pseudo-Riemannian calculus (\cite{aw:curvature.three.sphere}). Moreover, if $(C_{\A},h)$ is a free real metric calculus, then there exists a unique connection $\nabla$ such that $(C_{\A},h,\nabla)$ is pseudo-Riemannian (\cite{atn:nc.minimal.embeddings}).
\end{theorem}

Given the existence of a Levi-Civita connection, a noncommutative analogue to the classical Koszul formula can be used to find it in the case of a free real metric calculus. In the context of real connection calculi we state it as follows.

\begin{proposition}[\cite{aw:curvature.three.sphere}]
Let $(C_{\A},h,\nabla)$ be a pseudo-Riemannian calculus and assume that $\der_1,\der_2,\der_3\in\g$. Then it holds that
\begin{multline}\label{eqn:Koszul}
    2h(\nabla_{\der_1} e_2,e_3)=\der_1h(e_2,e_3)+\der_2h(e_1,e_3)-\der_3h(e_1,e_2)\\
    -h(e_1,\varphi([\der_2,\der_3]))+h(e_2,\varphi([\der_3,\der_1]))+h(e_3,\varphi([\der_1,\der_2])),
\end{multline}
where $e_i=\varphi(\der_i)$ for $i=1,2,3$.
Conversely, if $(C_{\A},h)$ is a real metric calculus and $\nabla$ is a connection satisfying Koszul's formula (\ref{eqn:Koszul}) for every $\der_1,\der_2,\der_3\in\g$, then $(C_{\A},h,\nabla)$ is a pseudo-Riemannian calculus.
\end{proposition}

\begin{example}
We shall briefly describe the construction of real calculi over the noncommutative torus $T^2_{\theta}$ and the noncommutative 3-sphere $S^3_{\theta}$. Since these are well-known examples of noncommutative manifolds it is interesting to see how they fit in the framework of real calculi. Both $T^2_{\theta}$ and $S^3_{\theta}$ are closely related to their respective classical counterparts $T^2$ and $S^3$ which are both parallelizable manifolds; this is reflected in the real calculi $C_{T^2_{\theta}}$ and $C_{S^3_{\theta}}$ constructed below being free. For more details about these constructions, see \cite{aw:curvature.three.sphere}.

The noncommutative torus $T^2_{\theta}$ is the $^*$-algebra with unitary generators $U,V$ satisfying the relation $VU=qUV$, where $q=e^{2\pi i\theta}$.
We let $\g'$ be the (real) Lie algebra generated by the hermitian derivations $\delta_1,\delta_2$, defined by
    \begin{align*}
        &\delta_1(U)=iU&\delta_2(U)&=0\\
        &\delta_1(V)=0&\delta_2(V)&=iV,
    \end{align*}
which implies that $[\delta_1,\delta_2]=0$.
In analogy with the classical torus $T^2$ being parallellizable, we let $M'$ be a free module of rank 2, with basis elements $e'_1,e'_2\in (T^2_{\theta})^4$.
With the map $\varphi':\delta_i\mapsto e'_i$ ($i=1,2$), we have the free real calculus $C_{T^2_{\theta}}=(T^2_{\theta},\g'_{D'},M',\varphi')$ over $T^2_{\theta}$ representing the noncommutative torus.
 
The noncommutative 3-sphere is the unital $^*$-algebra with generators $Z,Z^*,W,W^*$ subject to the relations
    \begin{align*}
      &WZ=qZW &W^*Z&=\bar{q}ZW^* &WZ^*&=\bar{q}Z^*W\\
      &W^*Z^*=qZ^*W^* &Z^*Z&=ZZ^* &W^*W&=WW^* \\
      &WW^*=\mathbbm{1}-ZZ^*.&&&&
    \end{align*}
We let $\g$ be the Lie algebra of dimension 3 generated by the hermitian derivations $\der_1,\der_2,\der_3$:
    \begin{align*}
      &\der_1(Z)=iZ, &\der_2(Z)&=0, &\der_3(Z)&=ZWW^*\\
      &\der_1(W)=0 &\der_2(W)&=iW &\der_3(W)&=-WZZ^*,
    \end{align*}
which implies that $[\der_i,\der_j]=0$ for every $i,j\in\{1,2,3\}$.
In analogy with the 3-sphere $S^3$ being parallellizable, we let $M$ be a free module of rank 3 with generators $e_1,e_2,e_3\in (S^3_{\theta})^4$ constituting a basis of $M$.
With the map $\varphi:\der_i\mapsto e_i$ ($i=1,2,3$), we have the free real calculus $C_{S^3_{\theta}}=(S^3_{\theta},\g_{D},M,\varphi)$ over $S^3_{\theta}$ representing the noncommutative 3-sphere.

Given invertible metrics $h':M'\times M'\rightarrow T^2_{\theta}$ and $h:M\times M\rightarrow S^3_{\theta}$ such that $(C_{T^2_{\theta}},h')$ and $(C_{S^3_{\theta}},h)$ are real metric calculi, Theorem~\ref{thm:existence.and.uniqueness.of.LC} guarantees that there are unique affine connections $\nabla$ and $\nabla'$ such that
$(C_{S^3_{\theta}},h,\nabla)$ and $(C_{T^2_{\theta}},h',\nabla')$ are pseudo-Riemannian. Since $C_{S^3_{\theta}}$ and $C_{T^2_{\theta}}$ are free, $\nabla$ and $\nabla'$ can be stated in terms of their respective Christoffel symbols which can be derived from Koszul's formula.
\end{example}

Given an algebra $\A$, a real Lie algebra $\g$ (isomorphic to a subalgebra of $\Der(\A)$ of hermitian derivations) and a (finitely generated) right $\A$-module $M$, it is interesting to ask how the choice of representation $D:\g\rightarrow\Der(\A)$ and linear map $\varphi:\g\rightarrow M$ affects the overall structure of the resulting real calculus $C_{\A}=(\A,\g_D,M,\varphi)$. To answer such questions a concept of real calculus homomorphisms is needed in order to determine when the structures of two real caluli are essentially the same, i.e., when they are isomorphic. 

\begin{definition}\label{def:rc.morphism}
	Let $C_{\A}$ and $C_{\A'}$ be real calculi and let $\phi:\A\rightarrow\A'$ be a $^*$-algebra homomorphism. The Lie algebra homomorphism $\psi:\g'\rightarrow\g$ is said to be compatible with $\phi$ if it satisfies the condition
	\begin{equation*}
	\phi(\psi(\delta)(a))=\delta(\phi(a))
	\end{equation*}
	for every $\der\in\g'$.
	If $\psi$ is compatible with $\phi$, then one sets $\Psi=\varphi\circ\psi$ and let $M_{\Psi}$ denote the submodule of $M$ generated by $\Psi(\g')$.
	Furthermore, if there exists a map $\psih:M_{\Psi}\rightarrow M'$ that satisfies the conditions
	\begin{enumerate}
		\item $\psih(m_1+m_2)=\psih(m_1)+\psih(m_2)$ for all $m_1,m_2\in M_{\Psi}$,\label{def:rc.morphism.psih.cond1}
		\item $\psih(ma)=\psih(m)\phi(a)$ for all $m\in M_{\Psi}, a\in\A$,\label{def:rc.morphism.psih.cond2}
		\item $\psih(\Psi(\der))=\varphi'(\der)$ for all $\der\in\g'$,\label{def:rc.morphism.psih.cond3}
	\end{enumerate}
	then $\psih$ is said to be compatible with $\phi$ and $\psi$, and we say that $(\phi,\psi,\psih)$ is a real calculus homomorphism from $C_{\A}$ to $C_{\A'}$; if $\phi$ and $\psi$ are isomorphisms and $\psih$ is a bijection, then $(\phi,\psi,\psih):C_{\A}\rightarrow C_{\A'}$ is called a real calculus isomorphism.
\end{definition}

Real calculus homomorphisms were originally developed in \cite{atn:nc.minimal.embeddings}, giving an important tool in the study of real calculi as algebraic structures. However, in the aforementioned article the main emphasis was placed on a geometric interpretation leading to a noncommutative theory of submanifolds. In what follows, the main focus is instead placed on algebraic aspects of real calculi, and we begin to explore the nuances involved in determining when two real calculi are essentially the same structure.

\section{Free and projective real metric calculi}\label{sec:free.and.projective}
    We shall consider real calculi $C_{\A}=(\A,\g_D,M,\varphi)$ where $M$ is projective, and we call such real calculi \emph{projective}. Given the importance of projective modules in algebra and geometry, it is interesting to know what general statements can be made about real calculi of this kind and in what way one should view them. Much like how every projective module can be realized as the projection of a free module, we show that a viable way to realize these structures is as "projections" of free real calculi. However, when we introduce metrics in an effort to construct a projective real metric calculus as a projection of a free real metric calculus the situation becomes more complicated, and additional conditions on the metric and projection in question must be satisfied.

In general, it is easy to generate projective real calculi from a given free real calculus.
\begin{proposition}\label{prop:derived.projective.calculus}
    Let $\tilde{C}_{\A}=(\A,\g_D,\A^n,\tilde{\varphi})$ be a free real calculus and let $P:\A^n\rightarrow \A^n$ be a projection. Then $C_{\A}=(\A,\g_D,P(\A^n),P\circ\tilde{\varphi})$ is a projective real calculus.
\end{proposition}
\begin{proof}
    All that needs to be checked is that the map $P\circ\tilde{\varphi}$ generates $P(\A^n)$ as a module over $\A$. But this is trivial, since $\tilde{\varphi}(\g)$ generates $\A^n$ and $P(m\cdot a)=P(m)\cdot a$, i.e., $P$ is a homomorphism. 
\end{proof}
In light of the above proposition, it is natural to ask whether every projective real calculus can be derived from a free real calculus in this way. This is indeed the case, as the next proposition shows.
\begin{proposition}\label{prop:projective.rc.realized.as.projection}
    Let $C_{\A}=(\A,\g_D,M,\varphi)$ be a projective real calculus where $\dim \g=n$. Then there exists a $\R$-linear map $\tilde{\varphi}:\g\rightarrow \A^n$   and a projection $P:\A^n\rightarrow\A^n$ such that
    \begin{equation*}
        (\A,\g_D,P(\A^n),P\circ\tilde{\varphi})\simeq C_{\A}.
    \end{equation*}
\end{proposition}
\begin{proof}
Let $\der_1,...,\der_n$ be a basis of $\g$ and let $e_k:=\varphi(\der_k)$ for $k=1,...,n$. Furthermore,
let $\hat{e}_1,...,\hat{e}_n$ be a basis of $\A^n$ and let $\tilde{\varphi}$ be the $\R$-linear map defined by $\tilde{\varphi}(\der_k):=\hat{e}_k$ for $k=1,...,n$. Define the module homomorphism $\rho:\A^n\rightarrow M$ by the formula $\rho(\hat{e}_kA^k)=e_kA^k$. Since $\varphi(\g)$ generates $M$, it follows that $\rho$ is an epimorphism. And since $M$ is projective it follows that $\rho$ splits, implying that there is a monomorphism $\nu:M\rightarrow \A^n$ such that $\rho\circ\nu=\text{id}_M$.
Let $P:=\nu\circ\rho$. Then it is clear that $P$ is a projection, since
\begin{equation*}
    P^2=(\nu\circ\rho)\circ(\nu\circ\rho)=\nu\circ(\rho\circ\nu)\circ\rho=\nu\circ\text{id}_M\circ\rho=\nu\circ\rho=P.
\end{equation*}
Let $\psih=P\circ \nu$. Then $\psih:M\rightarrow P(\A^n)$ is an isomorphism, with inverse $\psih^{-1}=\rho|_{P(\A^n)}$.

We now verify that $(\id{\A},\id{\g},\psih):C_{\A}\rightarrow (\A,\g_D,P(\A^n),P\circ\tilde{\varphi})$ is indeed a real calculus isomorphism. That $\text{id}_{\A}$ and $\text{id}_{\g}$ are compatible is trivial to check. Also, since $\psih$ is an isomorphism and since $\phi$ is the identity map, $\psih$ trivially satisfies conditions~(\ref{def:rc.morphism.psih.cond1}) and (\ref{def:rc.morphism.psih.cond2}) from Definition~\ref{def:rc.morphism}.

Finally, we have that
\begin{equation*}
    \psih(\varphi(\text{id}_{\g}(\der_k)))=\psih(e_k)=\psih(\rho(\hat{e}_k))=P\circ(\nu\circ\rho)(\hat{e}_k)=P^2(\hat{e}_k)=P\circ\tilde{\varphi}(\der_k)
\end{equation*}
for $k=1,...,n$, verifying that condition~(\ref{def:rc.morphism.psih.cond3}) for $\psih$ is satisfied as well. We conclude that $(\id{\A},\id{\g},\psih):C_{\A}\rightarrow (\A,\g_D,P(\A^n),P\circ\tilde{\varphi})$ is a real calculus isomorphism.
\end{proof}
\noindent
Thus, the above result says that any projective real calculus can be realized as a projection of a free real calculus $\tilde{C}_{\A}=(\A,\g_D,\A^n,\tilde{\varphi})$ in the sense that it is isomorphic to  $C_{\A}=(\A,\g_D,P(\A^n),P\circ\tilde{\varphi})$, where $P$ is a projection. 

Before discussing metric aspects of projective real calculi, we shall construct affine connections on these structures using the underlying structure of a free real calculus. As in differential geometry, there is a link between affine connections and \emph{bilinear maps}, i.e., maps $\Delta:\g\times M\rightarrow M$ that satisfy
\begin{align}
    &\Delta(\lambda\delta_1+\delta_2,m)=\lambda\Delta(\delta_1,m)+\Delta(\delta_2,m),\quad \delta_1,\delta_2\in\g, \lambda\in\R,\text{ and } m\in M,\label{eq:bilinear.map.1}\\ 
    &\Delta(\delta,m_1a+m_2)=\Delta(\delta,m_1)a+\Delta(\delta,m_2),\quad \delta\in\g, a\in\A,\text{ and } m_1,m_2\in M.\label{eq:bilinear.map.2}
\end{align}
The following lemma connecting bilinear maps and affine connections is a well-known result, and here we adapt it to our particular setting.
\begin{lemma}\label{lem:AffineConn.and.BilinForm}
    Let $\nabla^1,\nabla^2:\g\times M\rightarrow M$ be affine connections. Then $\Delta:=\nabla^1-\nabla^2$ is a bilinear map in the sense of (\ref{eq:bilinear.map.1}) and (\ref{eq:bilinear.map.2}). Converesely, if $\Delta:\g\times M\rightarrow M$ is a bilinear map and $\nabla:\g\times M\rightarrow M$ is an affine connection, then $\nabla':=\nabla+\Delta$ is an affine connection.
\end{lemma}
\begin{proof}
    For the first part, the only nontrivial condition to check is that $\Delta(\delta,ma)=\Delta(\delta,m)a$, but this readily follows from the Leibniz property:
    \begin{equation*}
        \Delta(\delta,ma)=\nabla^1_{\delta}(ma)-\nabla^2_{\delta}(ma)=(\nabla^1_{\delta}m-\nabla^2_{\delta}m)a+m\delta(a)-m\delta(a)=\Delta(\delta,m)a.
    \end{equation*}
    For the second part, the only nontrivial condition to check is the Leibniz condition. But this is a direct consequence of $\Delta$ being a bilinear form:
    \begin{align*}
        \nabla'_{\delta}(ma)&=\nabla_{\delta}(ma)+\Delta(\delta,ma)\\
        &=(\nabla_{\delta}m+\Delta(\delta,m))a+m\delta(a)=(\nabla'_{\delta}m)a+m\delta(a),
    \end{align*}
    and the proof is complete.
\end{proof}

We are now ready to state a standard result for connections on projective modules.
\begin{proposition}\label{prop:connection.on.projective.module}
    Let $P:\A^n\rightarrow\A^n$ be a projection. If $\tilde{\nabla}:\g\times\A^n\rightarrow \A^n$ is an affine connection, then the map $\nabla:\g\times P(\A^n)\rightarrow P(\A^n)$, defined by
    \begin{equation*}
        \nabla= P\circ\tilde{\nabla}
    \end{equation*}
    is an affine connection.
    Conversely, for every affine connection $\nabla:\g\times P(\A^n)\rightarrow P(\A^n)$ there is an affine connection $\tilde{\nabla}:\g\times\A^n\rightarrow \A^n$ such that $\nabla=P\circ \tilde{\nabla}$ on $\g\times P(\A^n)$.
\end{proposition}
\begin{proof}
	For the first statement one needs to check that $\nabla$ as defined is an affine connection, and since $P$ is a module homomorphism one only needs to check the Leibniz condition. Let $m\in P(\A^n)$, $\delta\in\g$ and $a\in \A$. Then
	\begin{equation*}
		\nabla_{\der} ma=P(\tilde{\nabla}_{\der}ma)=P((\tilde{\nabla}_{\der}m)a+m\der(a))=(\nabla_{\der}m)a+m\der(a)
	\end{equation*}
	since $P(m)=m$, which shows that $\nabla$ is indeed an affine connection.
	
    For the second statement, let $e_1,...,e_n$ be a basis of $\A^n$ and let $\tilde{\nabla}^0:\g\times \A^n\rightarrow \A^n$ be the flat connection with respect to this basis, i.e., $\tilde{\nabla}^0_{\der} (e_i a^i)=e_i\der(a^i)$. Next, let $\nabla^0:\g\times P(\A^n)\rightarrow P(\A^n)$ be the affine connection defined by
    \begin{equation*}
        \nabla^0_{\der}m=P(\tilde{\nabla}^0_{\der}m).
    \end{equation*}
    
    Let $\nabla:\g\times P(\A^n)\rightarrow P(\A^n)$ be an arbitrary affine connection. Then, by Lemma~\ref{lem:AffineConn.and.BilinForm}, it follows that $\Delta:=\nabla-\nabla^0$ is a bilinear map. Define the map $\tilde{\Delta}:\g\times\A^n\rightarrow\A^n$ by the formula
    \begin{equation*}
        \tilde{\Delta}(\der,m):=\Delta(\der,P(m));
    \end{equation*}
    then $\tilde{\Delta}$ is a bilinear map, since $P(ma)=P(m)a$.
    This, together with Lemma~\ref{lem:AffineConn.and.BilinForm}, implies that the map $\tilde{\nabla}:=\tilde{\nabla}^0+\tilde{\Delta}$ is an affine connection on $\g\times\A^n$. One quickly checks that $\nabla$ is indeed the projection of $\tilde{\nabla}$:
    \begin{align*}
        P(\tilde{\nabla}_{\der}m)&=P(\tilde{\nabla}^0_{\der}m+\tilde{\Delta}(\der,m))=\nabla^0_{\der}m+\Delta(\der,m)\\
        &=\nabla^0_{\der}m+(\nabla_{\der} m-\nabla^0_{\der}m)=\nabla_{\der}m,
    \end{align*}
    where the identity $P(\tilde{\Delta}(\der,m)))=\tilde{\Delta}(\der,m)=\Delta(\der,P(m))=\Delta(\der,m)$, $m\in P(\A^n)$, is used in the second equality above.
    This completes the proof.
\end{proof}

Let $\tilde{C}_{\A}=(\A,\g_D,\A^n,\tilde{\varphi})$ be a free real calculus and let $(\tilde{C}_{\A},\tilde{h})$ be a free real metric calculus. If $P:\A^n\rightarrow\A^n$ is a projection it has already been established in Proposition~\ref{prop:derived.projective.calculus} that we may generate the projective real calculus $C_{\A}=(\A,\g_D,P(\A^n),P\circ\tilde{\varphi})$ from $\tilde{C}_{\A}$. It is interesting to see under which conditions it is possible to generate a projective real metric calculus from $(\tilde{C}_{\A},\tilde{h})$, and if we assume that $P$ is orthogonal with respect to $\tilde{h}$ it is possible to make a few general statements.
\begin{proposition}\label{prop:restriction.of.metric}
	Let $\tilde{C}_{\A}=(\A,\g_D,\A^n,\tilde{\varphi})$ be a free real calculus, and let the projection $P$ be orthogonal with respect to the metric $\tilde{h}$ on $\A^n$, i.e., $\tilde{h}(P(m_1),m_2)=\tilde{h}(m_1,P(m_2))$ for all $m_1,m_2\in\A^n$. Then
	\begin{enumerate}
		\item The map $h(m_1,m_2)=\tilde{h}(m_1,m_2)$, $m_1,m_2\in P(\A^n)$, is a metric on $P(\A^n)$.
		\item If $\tilde{\nabla}:\g\times \A^n\rightarrow \A^n$ is an affine connection that is compatible with the metric $\tilde{h}$ then $\nabla=P\circ\tilde{\nabla}$ is compatible with $h$.  
	\end{enumerate}
\end{proposition}
\begin{proof}
	For the first statement, all that is needed is to show that $h$ is nondegenerate. Assume that $m_1\in P(\A^n)$ and that $h(m_1,m_2)=0$ for every $m_2\in P(\A^n)$. Then
	\begin{equation*}
		0=h(m_1,P(\tilde{m}))=\tilde{h}(m_1,P(\tilde{m}))=\tilde{h}(P(m_1),\tilde{m})=\tilde{h}(m_1,\tilde{m}),
	\end{equation*}
	where $\tilde{m}\in \A^n$ is arbitrary. Since $\tilde{h}$ is a metric it follows that $m_1=0$.
	
	For the second statement, Proposition~\ref{prop:connection.on.projective.module} ensures that $\nabla$ is an affine connection and if $\tilde{\nabla}$ is metric with respect to $\tilde{h}$ then a standard computation verifies that $\nabla$ is compatible with $h$ as well:
	\begin{align*}
		\der(h(m_1,m_2))&=\der(\tilde{h}(m_1,m_2))=\tilde{h}(\tilde{\nabla}_{\der} m_1,m_2)+\tilde{h}(m_1,\tilde{\nabla}_{\der}m_2)\\
		&=\tilde{h}(\tilde{\nabla}_{\der} m_1,P(m_2))+\tilde{h}(P(m_1),\tilde{\nabla}_{\der}m_2)\\
		&=\tilde{h}(\nabla_{\der} m_1,m_2)+\tilde{h}(m_1,\nabla_{\der}m_2)=h(\nabla_{\der} m_1,m_2)+h(m_1,\nabla_{\der}m_2),
	\end{align*}
	completing the proof.
\end{proof}

A natural question to ask is what conditions on $\tilde{h}$ and $P$ are needed to make sure that the projective real calculus
$C_{\A}=(\A,\g_D,P(\A^n),P\circ\tilde{\varphi})$ together with the restriction $h$ of $\tilde{h}$ to $P(\A^n)$ constitute a real metric calculus. Although Proposition~\ref{prop:restriction.of.metric} ensures that $h$ is a metric on $P(\A^n)$, it is not guaranteed that $h$ is symmetric on $\varphi(\g)$, i.e., $h(\varphi(\der_1),\varphi(\der_2))=h(\varphi(\der_2),\varphi(\der_1))$ for every $\der_1,\der_2\in\g$. Before presenting the next result we introduce some notation for the sake of clarity. Let $\{\der_1,...,\der_n\}$ be a basis of $\g$, and let $\hat{e}_i=\tilde{\varphi}(\der_i)$ and $e_i=(P\circ\tilde{\varphi})(\der_i)$ denote the generators of $\A^n$ and $P(\A^n)$, respectively. We introduce $\tilde{h}_{ij}$ and $h_{ij}$  to denote the components of the metric: 
\begin{equation*}
    \tilde{h}_{ij}=\tilde{h}(\hat{e}_i,\hat{e}_j),\quad h_{ij}=h(e_i,e_j).
\end{equation*}
Moreover, we introduce the coefficients $p^j_k\in\A$ as the unique elements such that
\begin{equation*}
    P(\hat{e}_k)=\hat{e}_j p^j_k;
\end{equation*}
since $P$ is a projection we have that $p^j_k=p^j_l p^l_k$.
\begin{proposition}
    Let $\tilde{h}$ be a metric on $\A^n$ such that $\tilde{h}_{ij}=\tilde{h}_{ji}$ with respect to the basis $\{\hat{e}_k\}_1^n$, and let $P$ be an orthogonal projection with respect to $\tilde{h}$. If $C_{\A}=(\A,\g_D,P(\A^n),\varphi)$ is a projective real calculus such that $\varphi(\der_k)=P(\hat{e}_k)=\hat{e}_l p^l_k$ and $h$ is the restriction of $\tilde{h}$ to $P(\A^n)$, then $(C_{\A},h)$ is a real metric calculus if and only if
    \begin{equation*}
        \tilde{h}_{jk}p^k_i=\tilde{h}_{ik}p^k_j
    \end{equation*}
    for any pair indices $i$ and $j$.
\end{proposition}
\begin{proof}
    From the orthogonality of $P$ it follows that
    \begin{equation}\label{eqn:orthogonal.metric.components}
        (p^k_i)^*\tilde{h}_{kj}=\tilde{h}(P(\hat{e}_i),\hat{e}_j)=\tilde{h}(\hat{e}_i,P(\hat{e}_j))=\tilde{h}_{ik}p^k_j.
    \end{equation}
    Using this, we see that
    \begin{equation*}
        h(\varphi(\der_i),\varphi(\der_j))=\tilde{h}(P(\hat{e}_i),P(\hat{e}_j))=(p^k_i)^*\tilde{h}_{kl}p^l_j\overset{(\ref{eqn:orthogonal.metric.components})}{=}\tilde{h}_{ik}p^k_lp^l_j,
    \end{equation*}
    and since $P$ is a projection (i.e., $P^2=P$), $\tilde{h}_{ik}p^k_lp^l_j$ can be simplified further to $\tilde{h}_{ik}p^k_j$.
    Thus $h(\varphi(\der_i),\varphi(\der_j))=h(\varphi(\der_j),\varphi(\der_i))$ is equivalent to $\tilde{h}_{ik}p^k_j=\tilde{h}_{jk}p^k_i$, and the statement follows.
\end{proof}
\section{Real calculi over matrix algebras}\label{sec:matrix.algebra}
    
Given two real calculi $C_{\A}=(\A,\g_D,M,\varphi)$ and $C_{\A'}=(\A',\g'_{D'},M',\varphi')$ it follows straight from the definition of real calculus isomorphisms that a necessary condition for $C_{\A}\simeq C_{\A'}$ is that $\A$ is isomorphic to $\A'$ and that $\g$ is isomorphic to $\g'$. Moreover, it was shown in \cite{atn:nc.minimal.embeddings} that the modules $M$ and $M'$ must also be isomorphic when viewed as modules over the same algebra (via the isomorphism between $\A$ and $\A'$).
These basic observations lead to the following natural question: given a unital $^*$-algebra $\A$, a Lie algebra $\g$ and a right $\A$-module $M$, how many nonisomorphic real calculi of the form $(\A,\g_D,M,\varphi)$ are there? Going forward, we shall refer to this as the classification problem for real calculi.

The classification problem for real calculi over general $^*$-algebras $\A$ is complicated, but the compatibility conditions between the maps $\phi$, $\psi$ and $\psih$ constituting an isomorphism $(\phi,\psi,\psih):(\A,\g_D,M,\varphi)\rightarrow (\A,\g_{D'},M,\varphi')$ imply that for a given $^*$-isomorphism $\phi:\A\rightarrow\A$ there is at most one choice of $\psi$ and $\psih$ such that $(\phi,\psi,\psih)$ is a real calculus homomorphism (cf. \cite{atn:nc.minimal.embeddings}); in general, the existence of a Lie algebra isomorphism $\psi:\g\rightarrow\g$ compatible with $\phi$ depends on the representations $D$ and $D'$, and the existence of a bijective map $\psih:M\rightarrow M$ that is compatible with $\phi$ and $\psi$ then depends on the relationship between the maps $\varphi$ and $\varphi'$. As a starting point we consider real calculi over matrix algebras, i.e., the case when $\A=\Mat{N}$. By the Skolem-Noether theorem all automorphisms of $\Mat{N}$ are inner, i.e., conjugations by a matrix $U\in\GL{N}$, and in the following we shall write $\phi_U$ to denote conjugation by $U\in\GL{N}$, i.e., $\phi_U(A)=U^{-1}AU$ for $A\in\Mat{N}$.

Another important property is that every derivation on $\Mat{N}$ is inner, i.e., any derivation $\der\in\text{Der}(\A)$ 
can be identified with the commutator of a unique trace-free matrix $\hat{D}\in\Mat{N}$:
\begin{equation*}\label{eqn:derivation.is.commutator.of.matrix}
    \der=[\hat{D},\cdot].
\end{equation*}
It follows that for any given representation $D:\g\rightarrow\text{Der}(\A)$ there is a unique matrix representation $\hat{D}:\g\rightarrow\Mat{N}$ such that $\hat{D}(\der)$ is trace-free and 
$D(\der)=[\hat{D}(\der),\cdot]$ for every $\der\in\g$. We call $\hat{D}$ the matrix representation of $\g$ associated with $D$, and note that since $D$ is faithful it follows that $\hat{D}$ is also a faithful representation of $\g$. Moreover, since every derivation $D(\delta)$ is assumed to be hermitian this has the effect that $\hat{D}(\delta)$ is anti-hermitian, implying that $\g$ is isomorphic to a Lie subalgebra of $\su{N}$. When discussing whether two real calculi $C_{\A}$ and $C'_{\A}$ are isomorphic when $\A=\Mat{N}$ 
one has to take into account that any real calculus isomorphism involves a Lie algebra isomorphism $\psi:\g\rightarrow\g$. This affects the relationship between the matrix representations $\hat{D}$ and $\hat{D}'$, and hence the concept of quasi-equivalent matrix representations becomes relevant. 

\begin{definition}
    Let $\hat{D}$ and $\hat{D}'$ be matrix representations of a Lie algebra $\g$. Then $\hat{D}$ and $\hat{D}'$ are said to be quasi-equivalent if there is a matrix $U\in\GL{N}$ and a Lie algebra automorphism $\psi\in\text{Aut}(\g)$ such that
    \begin{equation*}
        \hat{D}'(\der)=\phi_U(\hat{D}(\psi(\der)))
    \end{equation*}
    for every $\der\in\g$. The pair $(\phi_U,\psi)$ is called a \emph{realization} of the quasi-equivalence. 
\end{definition}
If the representations $\hat{D}$ and $\hat{D}'$ are faithful, then the automorphism $\psi$ in the above definition is unique since
\begin{align*}
    \phi_U\paraa{\hat{D}(\psi(\der))}=\phi_U\paraa{\hat{D}(\tilde{\psi}(\der))}\Leftrightarrow \hat{D}(\psi(\der))=\hat{D}(\tilde{\psi}(\der))\Leftrightarrow \psi(\der)=\tilde{\psi}(\der).
\end{align*}
Therefore, in this case we define $\psi_U$ to be the unique automorphism satisfying 
\begin{equation*}
    \hat{D}'(\der)=\phi_U(\hat{D}(\psi_U(\der)))
\end{equation*}
for every $\der\in\g$ whenever such an automorphism exists. In terms of real calculi over $\A=\Mat{N}$, quasi-equivalence of the matrix representations $\hat{D}$ and $\hat{D}'$ (associated with $D$ and $D'$ respectively) is the relevant condition for the existence of compatible automorphisms $\phi$ and $\psi$ (Definition~\ref{def:rc.morphism}), as can be seen in the following lemma.

\begin{lemma}\label{lem:Lie.automorphism.compatibility}
	Let $C_{\A}=(\Mat{N},\g_D,M,\varphi)$ and $C'_{\A}=(\Mat{N},\g_{D'},M,\varphi')$ be real calculi and suppose that $\phi\in\text{Aut}(\A)$ and $\psi\in\text{Aut}(\g)$. Then $\phi$ and $\psi$ are compatible (i.e., for every $a\in\A$ and $\der\in\g$, $\der(\phi(a))=\phi(\psi(\der)(a))$) if and only if $(\phi,\psi)$ is a realization of a quasi-equivalence between $\hat{D}$ and $\hat{D}'$.
\end{lemma}
\begin{proof}
	To prove sufficiency,
	assume that $\hat{D}$ and $\hat{D}'$ are quasi-equivalent.
	Then there is a nonsingular matrix $U$ such that $\hat{D}'(\der)=\phi_U(\hat{D}(\psi_U(\der)))$ for every $\der\in\g$. One readily checks that $\phi_U$ and $\psi_U$ are compatible:
	\begin{align*}
    	D'(\der)(\phi_U(A))&=[\hat{D}'(\der),\phi_U(A)]=[\phi_U(\hat{D}(\psi_U(\der))),\phi_U(A)]\\
    	&=\phi_U([\hat{D}(\psi_U(\der)),A])=\phi_U(D(\psi_U(\der))(A)).
    	\end{align*}
	To prove necessity, assume that $\phi:\Mat{N}\rightarrow \Mat{N}$ and $\psi:\g\rightarrow\g$ are compatible automorphisms. Then $\phi=\phi_U$ for a matrix $U\in\GL{N}$, and from the compatibility condition
	\begin{equation*}
	    \phi_U\paraa{D(\psi(\der))(A)}=D'(\der)(\phi_U(A))
	\end{equation*} 
	one obtains the following: 
	\begin{align*}
	[\phi_U(\hat{D}(\psi(\der))),\phi_U(A)]&=U^{-1}([\hat{D}(\psi(\der)),A])U=\phi_U(\psi(\der)(A))\\
	&=\der(\phi_U(A))=[\hat{D'}(\der),\phi_U(A)],\quad \der\in \g';
	\end{align*}
	since $\phi=\phi_U$ is an automorphism (and since $\hat{D}$ and $\hat{D}'$ are trace-free representations) this is equivalent to $\hat{D'}(\der)=\phi_U(\hat{D}(\psi(\der)))$ for all $\der\in \g$, implying that $\psi=\psi_U$, as desired.
\end{proof}
\noindent
In particular, the above lemma implies that if $C_{\A}$ and $C'_{\A}$ are isomorphic then the representations $\hat{D}$ and $\hat{D}'$ are quasi-equivalent.

Generally speaking, if only  free real calculi $C_{\A}=(\A,\g_D,M,\varphi)$ are considered where $\A,\g$ and $M$ are fixed and $D$ and $\varphi$ are allowed to vary, then the particular choice of $\varphi$ does not affect the overall structure of $C_{\A}$ in an essential way. This holds true not only for matrix algebras, but for general free real calculi. We state it as follows.
\begin{lemma}\label{lem:isomorphism.of.free.real.calculi}
	Let $C_{\A}=(\A,\g_D,M,\varphi)$ and $C'_{\A}=(\A,\g_{D'},M,\varphi')$ be free real calculi and suppose that $\phi\in\text{Aut}(\A)$ and $\psi\in\text{Aut}(\g)$ are compatible in the sense of Definition~\ref{def:rc.morphism}. Then $C_{\A}$ and $C_{\A'}$ are isomorphic as real calculi.
\end{lemma}
\begin{proof}
	Let $\{\delta_k\}_1^n$ be a basis of $\g$. Then $\{\der_k\}_1^n$, where $\der_k=\psi(\delta_k)$, is also a basis of $\g$. Moreover, since $C_{\A}$ and $C'_{\A}$ are free real calculi it follows that
	$\{\varphi(\der_k)\}_1^n$ and $\{\varphi'(\delta_k)\}_1^n$ are bases of $M$. 
	Let $\psih:M\rightarrow M$ and $\psit:M\rightarrow M$ be defined by the following:
	\begin{align*}
	&\psih(\varphi(\der_k)A^k):=\varphi'(\delta_k)\phi(A^k)\\
	&\psit(\varphi'(\delta_k)A^k):=\varphi(\der_k)\phi^{-1}(A^k).
	\end{align*}
	It is trivial to check that $(\phi,\psi,\psih):C_{\A}\rightarrow C'_{\A}$  and $(\phi^{-1},\psi^{-1},\psit):C'_{\A}\rightarrow C_{\A}$ are real calculus homomorphisms, and that $\psit=\psih^{-1}$.
\end{proof}
This result, together with Lemma \ref{lem:Lie.automorphism.compatibility}, can be used to describe how the structure of a free real calculus $C_{\A}=(\A,\g_D,M,\varphi)$ is completely determined by the representation $D$ (and its associated matrix representation $\hat{D}$) when $\A=\Mat{N}$ and $\g$ and $M$ are fixed.

\begin{theorem}\label{thm:free.real.calculi}
	Let 
	\begin{equation*}
	    C_{\A}=(\Mat{N},\g_D,M,\varphi)\quad\text{and}\quad C'_{\A}=(\Mat{N},\g_{D'},M,\varphi')
	\end{equation*}  
	be free real calculi. Then $C_{\A}\simeq C'_{\A}$ if and only if the matrix representations $\hat{D}$ and $\hat{D}'$ (associated with $D$ and $D'$, respectively) are quasi-equivalent.
\end{theorem}
\begin{proof}
	From Lemma \ref{lem:isomorphism.of.free.real.calculi} it follows that $C_{\A}\simeq C'_{\A}$ if and only if there are compatible automorphisms $\phi:\A\rightarrow\A$ and $\psi:\g\rightarrow\g$, and Lemma \ref{lem:Lie.automorphism.compatibility} states that such $\phi$ and $\psi$ exist if and only if $\hat{D}$ and $\hat{D}'$ are quasi-equivalent.
\end{proof}

Let $C_{\A}=(\Mat{N},\g_D,M,\varphi)$ and $C'_{\A}=(\Mat{N},\g_{D'},M,\varphi')$ be real calculi where $\text{dim }\g=n$.
Since Lemma \ref{lem:Lie.automorphism.compatibility} states that quasi-equivalence between the matrix representations $\hat{D}$ and $\hat{D}'$ is a necessary and sufficient condition for the existence of compatible automorphisms $\phi$ and $\psi$, this may be assumed when studying the role that the maps $\varphi$ and $\varphi'$ play in distinguishing the structures of $C_{\A}$ and $C'_{\A}$ from one another. In the case of free real calculi, Theorem \ref{thm:free.real.calculi} implies that the relationship between $\varphi$ and $\varphi'$ is irrelevant in order to determine whether two real calculi are isomorphic, but for general projective real calculi the situation changes. As an example we consider the case when $M=(\C^N)^n$, so that neither $C_{\A}$ nor $C'_{\A}$ is free whenever $n\neq N$.

Before moving forward we shall explain how the module $M=(\C^N)^n$ is represented for the sake of clarity. Since we are considering right modules, vectors in $\C^N$ are written as row vectors rather than column vectors, and $v\in (\C^N)^n$ is seen as a vector in $\C^{Nn}$ (in most cases written on the form $v=(v_1,...,v_n)\in(\C^N)^n$, with each $v_i\in\C^N$) , with matrices $A\in\Mat{N}$ acting on $(\C^N)^n$ in the following way:
\begin{equation*}
    v\cdot A=(v_1,...,v_n)\cdot A:=(v_1 A,...,v_n A)=v\left(\bigoplus_1^n A\right),
\end{equation*}
where $\oplus$ denotes the direct sum of matrices.
In general, the notation $v\cdot A$ above is used to distinguish the module action from a regular matrix multiplication when $n > 1$. Moreover, by $A\otimes B$ we will denote the Kronecker product of matrices.
\begin{proposition}\label{prop:projective.rc.iso}
    Let 
    \begin{equation*}
        C_{\A}=(\Mat{N},\g_D,(\C^N)^n,\varphi)\quad\text{and}\quad C'_{\A}=(\Mat{N},\g_{D'},(\C^N)^n,\varphi')
    \end{equation*} be real calculi. Then $C_{\A}$ and $C'_{\A}$ are isomorphic if and only if there are matrices $X\in\GL{n}$ and $U\in\GL{N}$ such that
    \begin{enumerate}
        \item $\hat{D}$ and $\hat{D}'$ are quasi-equivalent, with $\hat{D}'(\der)=\phi_U\paraa{\hat{D}(\psi_U(\der))}$,
        \label{prop:projective.rc.iso.cond1}
        \item $\varphi'(\der)=\varphi(\psi_U(\der))(X\otimes U)$\label{prop:projective.rc.iso.cond2}
    \end{enumerate}
    for all $\der\in\g$.
\end{proposition}
\begin{proof}
Assume that $(\phi,\psi,\psih):C_{\A}\rightarrow C'_{\A}$ is a real calculus isomorphism. Then Lemma \ref{lem:Lie.automorphism.compatibility} implies $\hat{D}$ and $\hat{D}'$ are quasi-equivalent and that there is a matrix $U\in\GL{N}$ such that $\phi=\phi_U$ and $\psi=\psi_U$.
By the linearity condition together with the compatibility condition with $\phi$, $\psih$  is a linear mapping from $(\C^N)^n$ to itself. Thus there is a unique matrix $\tilde{X}\in\Mat{Nn}$ such that $\psih(v)=v\tilde{X}$ for all $v\in (\C^N)^n$; for notational purposes, we identify $\tilde{X}$ with an $n$-by-$n$ matrix whose entries $\tilde{X}_{ij}$ are in $\Mat{N}$. Moreover, the compatibility condition between $\psih$ and $\phi$, $\psih(v\cdot A)=\psih(v)\cdot\phi_U(A)=(v\cdot A)\tilde{X}$, states that:
\begin{equation*}
    (v\cdot A)\tilde{X}=(v\tilde{X})\cdot (U^{-1}AU)\Leftrightarrow v\left[\tilde{X}\left(\bigoplus_1^n U^{-1}\right),\bigoplus_1^n A\right]=0
\end{equation*}
for all $v\in (\C^N)^n$ and $A\in\Mat{N}$. This is equivalent to $[\tilde{X}_{ij}U^{-1},A]=0$ for each $i,j,=1,...,n$, implying that each $\tilde{X}_{ij}=x_{ij}U$ for a constant $x_{ij}\in\C$; by setting $X:=(x_{ij})\in\Mat{n}$ we get the identity $\tilde{X}=X\otimes U$ and since $\psih$ is a bijection $X$ is invertible. Thus, $\psih:v\mapsto v(X\otimes U)$, and the compatibility condition between $\psih$ and $\psi_U$ states that
\begin{equation*}
    \varphi'(\der)=\psih(\varphi(\psi_U(\der)))=\varphi(\psi_U(\der))(X\otimes U)
\end{equation*}
for every $\der\in\g$, proving necessity of the stated conditions.

For sufficiency, we simply note that given matrices $X$ and $U$ satisfying the stated conditions we may define the bijective map $\psih:v\mapsto v(X\otimes U)$ on $(\C^N)^n$, and by the previous calculations it is clear that $\psih$ is compatible with both $\phi_U$ and $\psi_U$, implying that $(\phi_U,\psi_U,\psih)$ is a real calculus isomorphism from $(\Mat{N},\g_D,(\C^N)^n,\varphi)$ to $(\Mat{N},\g_{D'},(\C^N)^n,\varphi')$.
\end{proof}

A direct consequence of Proposition \ref{prop:projective.rc.iso} is that the classification task of real calculi $(\Mat{N},\g_D,(\C^N)^n,\varphi)$ where $D$ is given can be simplified in some cases by replacing $D$ with a quasi-equivalent representation $D'$ which is easier to work with. We state this as a lemma.

\begin{lemma}\label{lem:diagonal.form.of.matrix.representation}
Let $C_{\A}=(\Mat{N},\g_D,(\C^N)^n,\varphi)$ be a real calculus. Then, for any representation $D'$ that is  quasi-equivalent with $D$,
there is a $\R$-linear map $\varphi':\g\rightarrow (\C^N)^n$ 
such that $C_{\A}\simeq C'_{\A}=(\Mat{N},\g_{D'},(\C^N)^n,\varphi')$. 
\end{lemma}
\begin{proof}
    For any pair of quasi-equivalent representations $D'$ and  $D$ there is a matrix $U\in \GL{N}$ such that $\hat{D}'(\der)=\phi_U(\hat{D}(\psi_U(\der)))$, where $(\phi_U,\psi_U)$ is a realization of the quasi-equvialence.

    Choosing $\varphi':\g\rightarrow (\C^N)^n$ to be defined by the formula $\varphi'(\der)=\varphi(\psi_U(\der)))\cdot U$, it follows that $\varphi'$ is a linear map that generates $(\C^N)^n$ as a module over $\Mat{N}$. Moreover, since $\varphi'(\der)=\varphi(\psi_U(\der)))\cdot U=\varphi(\psi_U(\der)))(\One_n\otimes U)$, Proposition \ref{prop:projective.rc.iso} states that $C_{\A}\simeq C'_{\A}$, as desired.
\end{proof}

\subsection{The 1-dimensional case}\label{subsec:ex.dim.1}
When studying real calculi
\begin{equation*}
    C_{\A}=(\Mat{N},\g_{D},(\C^N)^n,\varphi)
\end{equation*} 
it is natural to begin with the simplest such case, i.e., when $n=\dim \g=1$. Although this may seem trivial at first glance, a more detailed survey of this case yields a lot of insight into what we can expect to find when the dimension of $\g$ is higher, especially in relation to metrics and connections. Also, due to the simplistic nature of $\operatorname{Aut}(\g)\simeq \Rstar$ in the one-dimensional case, it becomes possible to give a comprehensive discussion of the classification problem when the module is $\C^N$.

The analysis of isomorphism classes of real calculi $C_{\A}$ with a 1-dimensional Lie algebra $\g=\gen{\der}$ and a fixed representation $D:\g\rightarrow\Der(\Mat{N})$ is greatly simplified by the fact that all Lie algebra automorphisms of $\g$ are of the form $\psi_{\mu}(\der)=\mu\der$, where $\mu\in\Rstar$. Also, when $\g$ is 1-dimensional, it means that $\hat{D}$ is equivalent to a diagonal representation $\hat{D}'$, since $\hat{D}(\der)$ is skew-hermitian and thus diagonizable, with all eigenvalues being imaginary. 
Let
\begin{equation*}
    (\Mat{N},\gen{\der}_D,\C^N,\varphi)\quad\text{and}\quad (\Mat{N},\gen{\der}_D,\C^N,\varphi')
\end{equation*}
be real calculi. When analysing whether these are isomorphic, Proposition~\ref{prop:projective.rc.iso} implies that it is necessary to find every pair $U\in\GL{N}$ and $\mu\in\Rstar$ such that $\hat{D}(\der)=\phi_U(\hat{D}(\psi_{\mu}(\der)))=\mu\phi_U(\hat{D}(\der))$. To simplify this task, we may diagonalize $\hat{D}(\der)$ on the form
\begin{equation*}
    \hat{D}(\der)=\bigoplus_{j=1}^k \lambda_j I_{n_j},
\end{equation*}
where the eigenvalues $\lambda_j$ are sorted in a descending order (with respect to imaginary parts); by
Lemma~\ref{lem:diagonal.form.of.matrix.representation}, we may assume that $\hat{D}(\der)$ is of this diagonal form to begin with. In this particular case it is feasible to determine the possible pairs $U\in\GL{N}$ and $\mu\in\Rstar$ such that $\hat{D}(\der)=\mu\phi_U(\hat{D}(\der))$ directly. Since this implies that $\hat{D}(\der)$ and $\mu\hat{D}(\der)$ are similar, and since $\hat{D}(\der)$ is nonzero, this implies that $|\mu|=1$, with the special case $\mu=-1$ being possible if and only if $\hat{D}(\der)$ and $-\hat{D}(\der)$ are similar; in what follows we will refer to matrices $D$ such that $D$ and $-D$ are similar as \emph{anti-selfsimilar}. Furthermore, since $\hat{D}(\der)$ is diagonal with its eigenvalues sorted in a descending order, the condition $\hat{D}(\der)=\mu\phi_U(\hat{D}(\der))$ implies that $U$ is of the form
\begin{equation*}
    U_{+}=\begin{pmatrix}
      U_1 & &  \\
       & \ddots & \\
       &  & U_k
    \end{pmatrix}
    \quad\text{or}\quad
    U_{-}=\begin{pmatrix}
       & & U_1 \\
       & \udots & \\
      U_k &  & 
    \end{pmatrix},
\end{equation*}
where $U_j\in\GL{n_j}$ for $j=1,...,k$; we note in particular that the latter case where $U=U_{-}$ and $\mu=-1$ is possible only when $\hat{D}(\der)$ is anti-selfsimilar. 

\begin{lemma}\label{lem:diag.representation.isomorphism.criteria}
    Let     
    \begin{equation*}
        C_{\A}=(\Mat{N},\gen{\der}_D,\C^N,\varphi)\quad\text{and}\quad C'_{\A}=(\Mat{N},\gen{\der}_{D'},\C^N,\varphi')
    \end{equation*} 
    be real calculi such that
    \begin{equation*}
         \hat{D}(\der)=\bigoplus_{j=1}^k \lambda_j I_{n_j}\quad\text{and}\quad \hat{D}'(\der)=\mu_0\hat{D}(\der), \quad\mu_0\in\Rstar,
    \end{equation*}
    where the eigenvalues of $\hat{D}(\der)$ are sorted in a descending order (with respect to their imaginary parts)
    and where $n_j$ is the dimension of the eigenspace corresponding to the distinct eigenvalue $\lambda_j$ of $\hat{D}(\der)$. Moreover,
    write
    $\varphi(\der)=(v_1, v_2,..., v_k)$ and $\varphi'(\der)=(v'_1, v'_2,..., v'_k)$
    with $v_j,v'_j\in \C^{n_j}$ for $j=1,2,...,k$. If $\hat{D}(\der)$ is anti-selfsimilar, then 
    $C_{\A}$ is isomorphic to $C'_{\A}$ if and only if one of the two conditions hold:
    \begin{enumerate}
        \item $v_j=0\Leftrightarrow v'_j=0$ for every $j=1,2,...,k$,
        \item $v_j=0\Leftrightarrow v'_{k+1-j}=0$ for every $j=1,...,k$.
    \end{enumerate}
    If $\hat{D}(\der)$ is not anti-selfsimilar, then $C_{\A}$ is isomorphic to $C'_{\A}$ if and only if $v_j=0\Leftrightarrow v'_j=0$ for every $j=1,2,...,k$.
\end{lemma}
\begin{proof}
    Since $\hat{D}$ and $\hat{D}'$ are quasi-equivalent and $\gen{\der}$ is 1-dimensional, it is clear that every quasi-equivalence ($\phi_U$,$\psi_U$) between $\hat{D}$ and $\hat{D}'$ is such that $\psi_U:\der\mapsto \mu\der$. Thus, since $\hat{D}'(\der)=\mu_0\hat{D}(\der)=\mu\phi_U(\hat{D}(\der))$, it follows that $U\in\GL{N}$ satisfies $\phi_U(\hat{D}(\der))=\hat{D}(\der)$ (corresponding to $\mu=\mu_0$) or $\phi_U(\hat{D}(\der))=-\hat{D}(\der)$ (corresponding to $\mu=-\mu_0$); due to the specific diagonal form of $\hat{D}(\der)$, every such matrix $U$ can be written as
    \begin{equation*}
        U_{+}=\begin{pmatrix}
          U_1 & &  \\
           & \ddots & \\
           &  & U_k
        \end{pmatrix}
        \quad\text{or}\quad
        U_{-}=\begin{pmatrix}
           & & U_1 \\
           & \udots & \\
          U_k &  & 
        \end{pmatrix},
    \end{equation*}
    where $U_j\in\GL{n_j}$ for $j=1,...,k$.

    From this, together with Proposition \ref{prop:projective.rc.iso}, it follows that $C_{\A}\simeq C'_{\A}$ if and only if there is a matrix $U\in \GL{N}$ and $x\in\Rstar$ such that $v'=\mu x vU$ and such that $U$ is either of the form $U_{+}$ or $U_{-}$ as described above. This is equivalent to either
    \begin{enumerate}
        \item $v'_j=\mu_0 x v_j U_j$ for some $U_j\in \GL{n_j}$ for $j=1,2,...,k$, or
        \item $v'_j=-\mu_0 x v_{k+1-j} U_{k+1-j}$ for some $U_{k+1-j}\in \GL{n_{k+1-j}}$ for $j=1,2,...,k$;
    \end{enumerate}
    if $\hat{D}(\der)$ is anti-selfsimilar, then both of the above cases are possible, and if $\hat{D}(\der)$ is not, then only case (1) is possible. Since $U_j\in\GL{n_j}$ is arbitrary for $j=1,...,n$ in (1) and (2), the statement follows.
\end{proof}

Lemmas~\ref{lem:diagonal.form.of.matrix.representation} and~\ref{lem:diag.representation.isomorphism.criteria} imply the following corollary, which tells us the number of isomorphism classes there are for calculi of the form $(\Mat{N},\gen{\der}_D,\C^N,\varphi)$ when $D$ is given. In general, this number only depends on the eigenvalues of $\hat{D}(\der)$.
\begin{corollary}
Let $D:\g\rightarrow\Der(\Mat{N})$ be a faithful representation of $\g=\gen{\der}$, let $k$ be the number of distinct eigenvalues of $\hat{D}(\der)$ and let $|C_D|$ denote the number of pairwise nonisomorphic real calculi of the form $(\Mat{N},\g_D,\C^N,\varphi)$. Then
\begin{enumerate}
    \item if $\hat{D}(\der)$ is not anti-selfsimilar, then $|C_D|=2^k-1$,\\
    \item if $\hat{D}(\der)$ is anti-selfsimilar and $k$ is odd, then $|C_D|=2^{(k-1)/2}(1+2^{(k-1)/2})-1$,\\
    \item if $\hat{D}(\der)$ is anti-selfsimilar and $k$ is even, then $|C_D|=2^{k/2-1}(1+2^{k/2})-1$.
\end{enumerate}
\end{corollary}
\begin{proof}
    By Lemma~\ref{lem:diagonal.form.of.matrix.representation} we may assume that $\hat{D}(\der)$ is diagonal, i.e.,
    \begin{equation*}
        \hat{D}(\der)=\bigoplus_{j=1}^k \lambda_j I_{n_j}.
    \end{equation*}
    The only thing left to specify is $\varphi:\g\rightarrow M$, which must be chosen such that the image of $\g$ under $\varphi$ generates $M$; this means that the only restriction on $\varphi$ is that
    $\varphi(\der)=(v_1, v_2,...,v_k)\neq 0$. This restriction is the same as requiring that at least one of $v_1,v_2,...,v_k$ has to be chosen to be nonzero. For $v=(v_1, v_2,...,v_k)\in\C^N$, let $a(v)=(a_1,a_2,...,a_k)\in \{0,1\}^k$ be such that $a_j=1$ if $v_j\neq 0$ and $a_j=0$ if $v_j=0$, and set $\tilde{a}(v)=(\tilde{a}_1,...,\tilde{a}_k)=(a_k,...,a_1)$. 
    
    If $\hat{D}(\der)$ is not anti-selfsimilar, Lemma \ref{lem:diag.representation.isomorphism.criteria} states that 
    the choices $\varphi(\der)=v$ and $\varphi'(\der)=v'$
    yield isomorphic real calculi if and only if
    $a(v')=a(v)\in\{0,1\}^k$, and since there are exactly $2^k-1$ nonzero elements in $\{0,1\}^k$, $|C_D|=2^k-1$ in this case.
    
    If $\hat{D}(\der)$ is anti-selfsimilar one gets two relevant cases.
    Suppose that $k=2m+1$. Then there are exactly $2^{m+1}$ elements $b\in\{0,1\}^k$ such that $b=(b_1,...,b_k)=(b_k,...,b_1)$, implying that there are exactly $2^k-2^{m+1}$ elements $c\in\{0,1\}^k$ such that $c=(c_1,...,c_k)\neq (c_k,...,c_1)$. By Lemma \ref{lem:diag.representation.isomorphism.criteria},
    the choices $\varphi(\der)=v$ and $\varphi'(\der)=v'$
    yield isomorphic real calculi if and only if either
    $a(v')=a(v)\in\{0,1\}^k$ or $a(v')=\tilde{a}(v)\in\{0,1\}^k$; this distinction does not matter if $a(v)=\tilde{a}(v)$, and as a result from the above we have
    that the number of pairwise nonisomorphic real calculi of the given form is
    \begin{equation*}
        \frac{2^k-2^{m+1}}{2}+(2^{m+1}-1)=2^2m-2^m+2\cdot 2^m-1=2^m(2^m+1)-1,
    \end{equation*}
    since $a(v)=(0,...,0)$ does not correspond to a valid choice of $v\in\C^N$.
    
    If, instead, $k=2m$, then there are exactly $2^m$ elements $b\in\{0,1\}^k$ such that $b=(b_1,...,b_k)=(b_k,...,b_1)$, implying that there are exactly $2^k-2^m$ elements $c\in\{0,1\}^k$ such that $c=(c_1,...,c_k)\neq (c_k,...,c_1)$. Thus, the number of pairwise nonisomorphic real calculi of the given form is
    \begin{equation*}
        \frac{2^k-2^m}{2}+(2^m-1)=2^{2m-1}-2^{m1}+2\cdot 2^{m-1}-1=2^{m-1}(2^m+1)-1
    \end{equation*}
    in this case, completing the proof.
\end{proof}

Next, we wish to study real metric calculi and connections and in the spirit of Section~\ref{sec:free.and.projective} we shall do this by realizing the real calculus $C_{\A}=(\Mat{N},\g_D,\C^N,\varphi)$ as a projection of a free real calculus. By Proposition~ \ref{prop:projective.rc.realized.as.projection}
this amounts to finding a projection $P:\A\rightarrow\A$ and a map $\tilde{\varphi}:\g\rightarrow\A$ such that $P(\A)\simeq \C^N$ and such that $\varphi$ can be identified with $P\circ\tilde{\varphi}$. If we let $P:\A\rightarrow\A$ be the projection that maps a matrix $A$ to the matrix $A'$ whose first row is equal to that of $A$ and whose other rows are zero, and let $\hat{e}$ be an invertible matrix whose first row is $v_0=\varphi(\der)$, we may set $\tilde{\varphi}(\der)=\hat{e}$; by construction, these choices of projection $P$ and map $\tilde{\varphi}$ yield us a free real calculus $\tilde{C}_{\A}$ from which we may realize $C_{\A}$ as a projection. For the sake of simplicity, $\hat{e}$ is chosen so that its rows are pairwise orthogonal, and for the sake of clarity we make the following identification between $\C^N$ and $P\paraa{\Mat{N}}$:
\begin{equation*}
    (x_1,...,x_N)\simeq \theta\paraa{(x_1,...,x_N)}:=\begin{pmatrix}
      x_1 & \cdots & x_N \\
      0 & \cdots & 0\\
      \vdots  & \ddots & \vdots \\
      0 & \cdots & 0
    \end{pmatrix},
\end{equation*}
where the map $\theta:\C^N\rightarrow P(\Mat{N})$ defined above is an obvious isomorphism.
Although the above description of $P$ may seem clear, since $P$ is going to be used with respect to the basis of $\A$ given by $\hat{e}$ rather than the standard basis, a more detailed description of $P$ is warranted. When viewed in the basis $\{\hat{e}\}$, $P$ can be regarded as the module homomorphism satisfying
\begin{equation*}
    P(\hat{e})=\hat{e}p,
\end{equation*}
where $p=\frac{1}{||v_0||^2}v_0^{\hconj}v_0$. Since the rows of $\hat{e}$ are chosen to be pairwise orthogonal and since $p=p^2$ is the matrix that projects vectors in $\C^N$ onto $v_0$, it is clear that $P$ is indeed the projection described earlier, with $\varphi(\der)=v_0\simeq\theta(v_0)=P(\hat{e})=(P\circ\tilde{\varphi})(\der)$.

Next, let us show that every metric on $\C^N$ is of the form $h(u,v)=x\cdot u^{\hconj}v\in\Mat{N}$.

\begin{proposition}\label{prop:metrics.on.CN}
    Let $h$ be a metric on the right $\Mat{N}$-module $\C^N$. Then there exists $x\in\Rstar$ such that 
    \begin{equation*}
        h_x(u,v)= x \cdot u^\hconj v\in\Mat{N}.
    \end{equation*}
    for $u,v\in\C^N$.
\end{proposition}
\begin{proof}
Let $h$ be a metric on $\C^N$, and let $e_1=(1,0,...,0)\in\C^N$. Then we may set $H:=h(e_1,e_1)$ and use the identity $v=e_1\theta(v)$ (where $\theta$ is the isomorphism between $\C^N$ and $P(\Mat{N})$ described earlier) to calculate $h(u,v)$ for general vectors in $\C^N$:
\begin{align*}
    h(u,v)&=h(e_1\theta(u),e_1\theta(v)=\theta(u)^{\hconj} h(e_1,e_1)\theta(v)=\theta(u)^{\hconj} H\theta(v)\\
    &=\begin{pmatrix}
      \overline{u}_1  & 0 &\cdots & 0\\
      \overline{u}_2  & 0 &\cdots & 0\\
      \vdots  & \vdots & \ddots & \vdots\\
      \overline{u}_N  & 0 &\cdots & 0\\
    \end{pmatrix}
    \begin{pmatrix}
      h_{11}  & h_{12} &\cdots & h_{1N}\\
      h_{21}  & h_{22} &\cdots & h_{2N}\\
      \vdots  & \vdots & \ddots & \vdots\\
      h_{N1}  & h_{N2} &\cdots & h_{NN}\\
    \end{pmatrix}
    \begin{pmatrix}
      v_1  & v_2 &\cdots & v_N\\
      0  & 0 &\cdots & 0\\
      \vdots  & \vdots & \ddots & \vdots\\
      0  & 0 &\cdots & 0\\
    \end{pmatrix}\\
    &=\begin{pmatrix}
      \overline{u}_1 \\
      \overline{u}_2\\
      \vdots \\
      \overline{u}_N \\
    \end{pmatrix}
     (h_{11}, h_{12},\cdots,h_{1N})
    \begin{pmatrix}
      v_1  & v_2 &\cdots & v_N\\
      0  & 0 &\cdots & 0\\
      \vdots  & \vdots & \ddots & \vdots\\
      0  & 0 &\cdots & 0\\
    \end{pmatrix}=h_{11}\cdot u^{\hconj}v.
\end{align*}
If $h_{11}$ is zero, then $h(u,v)=0$ for every $u,v\in\C^N$. Thus, $h_{11}\in\Rstar$ and the statement follows.
\end{proof}

We now move on to affine connections $\nabla$ on $\C^N$; by Proposition \ref{prop:connection.on.projective.module} every such connection is the composition of the projection $P$ (as defined earlier with respect to the basis $\hat{e}=\tilde{\varphi}(\der)$ of $\Mat{N}$) with a connection $\tilde{\nabla}$ on the free module $\A$ and thus we may use the structure of the free calculus $\tilde{C}_{\A}$ to define $\nabla$. On the free module it is straightforward to define an affine connection in terms of its Christoffel symbol $\tilde{\Gamma}\in\A$ (with respect to the basis $\hat{e}=\tilde{\varphi}(\der)$):
\begin{equation*}
    \tilde{\nabla}_{\der} \hat{e}=\hat{e}\tilde{\Gamma},
\end{equation*}
and for the connection $\nabla=P\circ\tilde{\nabla}$ on $P(\Mat{N})$ we have
\begin{equation*}
    \nabla_{\der} P(\hat{e})=P(\tilde{\nabla}_{\der} P(\hat{e}))=P(\tilde{\nabla}_{\der} \hat{e}p)=P(\hat{e}\tilde{\Gamma})p+P(\hat{e})\der(p)=P(\hat{e})(\tilde{\Gamma}p+\der(p)),
\end{equation*}
where we recall that
\begin{equation*}
    p=\frac{1}{||v_0||^2}v_0^{\hconj}v_0
\end{equation*}
is the matrix projecting vectors in $\C^N$ onto $v_0=\varphi(\der)$. Thus, through the identification $v_0\simeq \theta(v_0)=P(\hat{e})$, we see that any connection $\nabla$ on $\C^N$ is defined by
\begin{equation}\label{eqn:connection.on.CN.first.form}
    \nabla_{\der} v_0=v_0(\tilde{\Gamma} p+\der(p)),
\end{equation}
where $\tilde{\Gamma}\in\Mat{N}$ can be chosen arbitrarily.
This can be simplified somewhat since the product
\begin{equation*}
    v_0\tilde{\Gamma} p=\left(\frac{1}{||v_0||^2}v_0\tilde{\Gamma} v_0^{\hconj}\right)v_0
\end{equation*}
is equal to $\lambda v_0$ for a unique $\lambda\in\C$ for every $\tilde{\Gamma}\in\A$, implying that we may define $\nabla$
\begin{equation}\label{eqn:connection.on.CN}
    \nabla_{\der} v_0=v_0(\lambda\One_N+\der(p))
\end{equation}
for an arbitrary $\lambda\in\C$. Note that an arbitrary connection could also be characterized by 
\begin{equation*}
    \nabla_{\der} v_0=v_0(\tilde{\lambda}\One_N-\hat{D}(\der))
\end{equation*}
for an arbitrary $\lambda\in\C$, since equation~(\ref{eqn:connection.on.CN.first.form}) is equivalent to
\begin{align*}
    \nabla_{\der} v_0&=v_0(\tilde{\Gamma} p+\der(p))=v_0(\tilde{\Gamma} p+[\hat{D}(\der),p])=v_0((\tilde{\Gamma}+\hat{D}(\der))p-p\hat{D}(\der))\\
    &=v_0(\tilde{\lambda}\One-\hat{D}(\der)).
\end{align*}

At this point we have managed to define a connection on the projective module, and
by Proposition \ref{prop:metrics.on.CN} we know that any metric $h$ on $\C^N$is given by $h=h_x:(u,v)\mapsto x\cdot u^{\hconj}v$. One easily checks that the pair $(C_{\A},h_x)$ forms a real metric calculus for any real $x\neq 0$, and so we may use the above characterization of connections on $\C^N$ to determine the possible real connection calculi on the form $(C_{\A},h_x,\nabla)$. In general, the existence of a real connection calculus $(C_{\A},h_x,\nabla)$ relies upon $v_0=\varphi(\der)$ being an eigenvector of $\hat{D}(\der)$, as can be seen in the following proposition.

\begin{proposition}\label{prop:v0.is.eigenvector.and.character.of.nabla}
    Let $\hat{D}=\hat{D}(\der)$ be the matrix representation of $\der\in\g$. Then, for any $x\in\Rstar$, there exists a connection $\nabla$ such that
    $(C_{\A},h_x,\nabla)$ is a real connection calculus if and only if $v_0$ is an eigenvector of $\hat{D}(\der)$. In this case, if $\nabla^0$ is a given connection on $\C^N$, then $(C_{\A},h_x,\nabla^0)$ is a real connection calculus if and only if $\nabla^0_{\der}v_0=\lambda v_0$ for $\lambda\in\R$.
\end{proposition}

\begin{proof}
Assume that $(C_{\A},h_x,\nabla)$ is a real connection calculus.
Using equation~(\ref{eqn:connection.on.CN}), we have that $\nabla_{\der}v_0=v_0(\lambda\One_N+\der(p))$ for a $\lambda\in\C$. Thus, we find that the symmetry condition for real connection calculi, $h_x(v_0,\nabla_{\der}v_0)=h_x(v_0,\nabla_{\der}v_0)^{\hconj}$, is equivalent to
\begin{equation*}
    h_x(v_0,\nabla_{\der}v_0)=h_x(v_0,v_0(\lambda\One_N+\der(p)))=x\cdot v_0^{\hconj}v_0(\lambda\One_N+\der(p))=x\cdot(\bar{\lambda}\One_N+\der(p))v_0^{\hconj}v_0,
\end{equation*}
since $\der$ is a hermitian derivation and $p=p^{\hconj}$. Using $v_0^{\hconj}v_0=||v_0||^2 p$, one gets
\begin{align*}
    (\lambda-\bar{\lambda})v_0^{\hconj}v_0&=\der(p)v_0^{\hconj}v_0-v_0^{\hconj}v_0\der(p)\\
    &=||v_0||^2([\hat{D},p]p-p[\hat{D},p])=||v_0||^2(\hat{D}p+p\hat{D}-2p\hat{D}p).
\end{align*}
Multiplying the above expression by $v_0$ from the left, we note that $v_0$ is a (left) eigenvector of $\hat{D}p+p\hat{D}-2p\hat{D}p$ with eigenvalue $(\lambda-\bar{\lambda})$, meaning that
\begin{align*}
    (\lambda-\bar{\lambda})v_0&=v_0(\hat{D}p+p\hat{D}-2p\hat{D}p)=v_0\hat{D}(\One_N-p),
\end{align*}
since $v_0 p=v_0$.
Multiplying this expression from the right by $v_0^{\hconj}$, we find that
\begin{equation*}
    (\lambda-\bar{\lambda})||v_0||^2=v_0\hat{D}(v_0^{\hconj}-v_0^{\hconj})=0,
\end{equation*}
which in turn implies that
\begin{equation*}
    v_0\hat{D}(\One_N-p)=0.
\end{equation*}
Hence, $v_0\hat{D}$ lies in the kernel of $(\One_N-p)$ which is the same as saying that $v_0$ is an eigenvector of $\hat{D}$.

To complete the proof, let $v_0$ be an eigenvector of $\hat{D}$ with eigenvalue $i \mu$ (where $\mu\in\R$ since $\hat{D}$ is anti-hermitian), and let $\nabla^0_{\der} v_0=v_0(\lambda\One_N+\der(p))$, where $\lambda\in\C$. Since $p=\frac{1}{||v_0||^2}v_0^{\hconj}v_0$ and $\hat{D}^{\hconj}=-\hat{D}$, we have that
\begin{equation*}
    \hat{D}p=-\frac{1}{||v_0||^2}(v_0\hat{D})^{\hconj}v_0=-\frac{1}{||v_0||^2}(i\mu v_0)^{\hconj}v_0 =i\mu p,  
\end{equation*}
which implies that $\der(p)=\hat{D}p-p\hat{D}=i\mu p-i\mu p=0$. Thus, $\nabla^0_{\der} v_0=\lambda v_0$, and the symmetry condition $h_x(v_0,\nabla^0_{\der}v_0)=h_x(v_0,\nabla^0_{\der}v_0)^{\hconj}$ then becomes equivalent to
\begin{equation*}
   x\lambda\cdot v_0^{\hconj}v_0=x\bar{\lambda}\cdot v_0^{\hconj}v_0.
\end{equation*}
This is true if and only if $\lambda\in\R$.
\end{proof}

We are now ready to state a general result regarding the existence of a Levi-Civita connection given a real metric calculus $(C_{\A},h_x)$.

\begin{proposition}
    Let $C_{\A}=(\Mat{N},\gen{\der}_D,\C^N,\varphi)$ and let $(C_{\A},h_x)$ be a real metric calculus.
    Then there exists a unique connection $\nabla$ such that $(C_{\A},h_x,\nabla)$ is pseudo-Riemannian if and only if $v_0=\varphi(\der)$ is an eigenvector of $\hat{D}(\der)$ with eigenvalue $\lambda_{\der}$. In this case, the Levi-Civita connection $\nabla$ is given by
    \begin{equation*}
        \nabla_{\der} v=\lambda_{\der}v-v\hat{D}(\der),
    \end{equation*}
    for $v\in\C^N$.
\end{proposition}
\begin{proof}
Assume that $(C_{\A},h_x,\nabla)$ is a real connection calculus.
Then Proposition~\ref{prop:v0.is.eigenvector.and.character.of.nabla} immediately implies that $v_0$ is an eigenvector of $\hat{D}$ with eigenvalue $\lambda_{\der}$. Using the fact that $\nabla_{\der} v_0=\lambda v_0$ for some $\lambda\in\R$, and noting that $\der(h_x(v_0,v_0))=0$ (since $h_x(v_0,v_0)$ is proportional to $p$, and since $\der(p)=0$), the metric condition for $\nabla$ becomes
\begin{align*}
    0=\der(h_x(v_0,v_0))=h_x(v_0,\nabla_{\der}v_0)+h_x(\nabla_{\der}v_0,v_0)=(\lambda+\bar{\lambda}) h_x(v_0,v_0).
\end{align*}
This is satisfied if and only if $\lambda+\bar{\lambda}=0$, and since $\lambda\in\R$ it follows that $2\lambda=0$. In the case of a one-dimensional Lie algebra the torsion condition is trivial, which means that the real connection calculus $(C_{\A},h_x,\nabla)$ is pseudo-Riemannian if and only if $\nabla$ is given by
\begin{equation*}
    \nabla_{\der} v=\nabla_{\der} (v_0 B) =v_0\der(B)=v_0[\hat{D}(\der),B]=\lambda_{\der} v_0 B-(v_0 B)\hat{D}(\der)=\lambda_{\der} v-v\hat{D}(\der),
\end{equation*}
where $v=v_0B$ is an arbitrary vector in $\C^N$.
\end{proof}

\subsection{General abelian Lie algebras}\label{subsec:ex.dim.n}
Let $C_{\A}=(\Mat{N},\g_D,(\C^N)^n,\varphi)$ and assume that $\g$ is abelian with a basis $\{\der_1,...,\der_n\}$. Some of the results from the previous section can be generalized to higher dimensions without great effort, and to this end we will consider the special case where each $e_i:=\varphi(\der_i)$ is of the following form:
\begin{equation*}
    e_i=(0,...,0,\alpha_i v_0,0,...,0)\in(\C^{N})^n,
\end{equation*}
where $v_0\in\C^{N}$ is a unit vector and $\alpha_i\in\Rstar$. The effect this has is that for any metric $h$ on $(\C^{N})^n$ we have that
$h(e_i,e_j)=\mu_{ij}\alpha_i\alpha_j v_0^\hconj v_0=\tilde{\mu}_{ij}v_0^\hconj v_0$ where $\mu_{ij}\in\R$ (implying that $\tilde{\mu}_{ij}\in\R$ as well); this can be most easily seen by noting that $h$ can be broken down as the sum of general bilinear forms over $\C^N$, after which we can apply the exact same reasoning used to prove Proposition \ref{prop:metrics.on.CN} to reach this conclusion. Moreover, the symmetry condition $h(\varphi(\der_i),\varphi(\der_j))=h(\varphi(\der_j),\varphi(\der_i))$
implies that $\tilde{\mu}_{ij}=\tilde{\mu}_{ji}\in\R$ for every pair $i,j$ and the metric $h$ being nondegenerate is reflected in the fact that the real and symmetric $n\times n$ matrix $\tilde{M}=(\tilde{\mu}_{ij})$ defining the metric $h$ is invertible.

To realize $(C_{\A},h)$ as a projection of a free real metric calculus we need to determine a suitable projection $P:\A^n\rightarrow\A^n$ and invertible metric $\tilde{h}$ on $\A^n$, and to this end we let $A$ be a unitary matrix whose first row is equal to $v_0$. We then set $A_i=\alpha_i A$ and let $\{\hat{e}_i\}_1^n$ be the basis of $\A^n$ given by
\begin{equation*}
    \hat{e}_i=(0,...,0,A_i,0,...,0)\in\A^n.
\end{equation*}
As a projection we use the module homomorphism $P:\A^n\rightarrow\A^n$ given by
\begin{equation*}
    P(\hat{e}_j)=\hat{e}_kp^k_j,
\end{equation*}
where
\begin{equation*}
    p^k_j=\begin{cases}
    v_0^{\hconj}v_0=p &\text{if }j=k,\\
    0 &\text{otherwise}.
    \end{cases}
\end{equation*}
Since $p$ is the matrix that projects vectors in $\C^N$ onto $v_0$, and since each $A_i=\alpha_i A$ where $A$ is a unitary matrix whose first row is $v_0$, the product $A_i p$ is the matrix whose first row is $\alpha_i v_0$ and whose other rows are all zero.
Using the map $\theta:\C^N\rightarrow\Mat{N}$ from the previous section, given by
\begin{equation*}
     \theta\paraa{(x_1,...,x_N)}:=\begin{pmatrix}
      x_1 & \cdots & x_N \\
      0 & \cdots & 0\\
      \vdots  & \ddots & \vdots \\
      0 & \cdots & 0
    \end{pmatrix},
\end{equation*}
we define the map $\Theta:(\C^N)^n\rightarrow P(\A^n)$ as
\begin{equation*}
    \Theta(v_1,...,v_n)=(\theta(v_1),...,\theta(v_n)),
\end{equation*}
where each $v_i\in\C^N$. It is easy to check that $\Theta$ is an isomorphism between $(\C^N)^n$ and $P(\A^n)$, so that vectors $v\in (\C^N)^n$ can be identified with $\Theta(v)\in P(\A^n)$. Using this identification we note that $\varphi(\der_i)=e_i\simeq\Theta(e_i)=P(\hat{e}_i)$ for $i=1,...,n$, and thus we define the map $\tilde{\varphi}:\g\rightarrow \A^n$ by $\tilde{\varphi}(\der_i)=\hat{e}_i$ for $i=1,...,n$. This yields the free real calculus $\tilde{C}_{\A}=(\Mat{N},\g_D,(\Mat{N})^n,\tilde{\varphi})$ from which $C_{\A}$ can be realized as a projection.

Next, we choose the metric $\tilde{h}$ to be given by
\begin{equation*}
    \tilde{h}(\hat{e}_i,\hat{e}_j)=\tilde{h}_{ij}:=\tilde{\mu}_{ij}A^{\hconj}A=\tilde{\mu}_{ij}\One_N.
\end{equation*}
One quickly verifies that this is an invertible metric by calculating $\tilde{h}^{ij}\tilde{h}_{jk}$, where $\tilde{h}^{ij}=\tilde{\mu}^{ij}\One_N$. It is straightforward to verify that $P$ is orthogonal with respect to $\tilde{h}$ and that $h(e_i,e_j)=\tilde{h}(P(\hat{e}_i),P(\hat{e}_j))$, which means that $(C_{\A},h)$ can be realized as the projection of the free real metric calculus $(\tilde{C}_{\A},\tilde{h})$.

We are interested in affine connections $\nabla$ on $(\C^{N})^n$, and since every such connection is the projection of a connection $\tilde{\nabla}$ on the free module, we may use the structure on the free module to perform calculations.
On the free module, it is trivial to define an affine connection $\tilde{\nabla}$ in terms of its Christoffel symbols $\tilde{\Gamma}^k_{ij}$, where 
\begin{equation*}
    \nabla_i \hat{e}_j=\hat{e}_k\tilde{\Gamma}^k_{ij}.
\end{equation*}
Using this, we may calculate corresponding Christoffel symbols $\Gamma^k_{ij}$ for the connection $\nabla:=P\circ\tilde{\nabla}$:
\begin{align*}
    \nabla_i e_j&=P(\tilde{\nabla}_i P(\hat{e}_j))=P(\tilde{\nabla}_i\hat{e}_k p^k_j)=P(\tilde{\nabla}_i\hat{e}_j p)=P(\hat{e}_k\tilde{\Gamma}^k_{ij}p+\hat{e}_j\der_i(p))\\
    &=e_k\tilde{\Gamma}^k_{ij}p+e_j\der_i(p)=e_k(\tilde{\Gamma}^k_{ij}p+\delta^k_j\der_i(p)),
\end{align*}
where $p=v_0^{\hconj}v_0$ is the matrix projecting vectors in $\C^N$ onto $v_0$.
Thus, if we define $\nabla$ by
\begin{equation*}
    \nabla_i e_j:=e_k\Gamma^k_{ij},
\end{equation*}
where each $\Gamma^k_{ij}=\tilde{\Gamma}^k_{ij}p+\delta^k_j\der_i(p)$ for a $\tilde{\Gamma}^k_{ij}\in\A$, then $\nabla$ is a well-defined connection on $(\C^N)^n$. We are, however, only interested in connections such that $(C_{\A},h,\nabla)$ is a real connection calculus, and for this to be the case it is necessary that we assume that $v_0$ is an eigenvector of every $\hat{D}(\der_i)$; this requirement is proven in the same way that Proposition~\ref{prop:v0.is.eigenvector.and.character.of.nabla} is proved in the 1-dimensional case. Since we are only considering abelian Lie-algebras $\g$ such a vector $\tilde{v}_0$ of course always exists, and as in the 1-dimensional case this has the consequence that $\der_i(p)=\der_i(h_{kj})=0$. Moreover, we note that for every matrix $\tilde{\Gamma}^k_{ij}$ there is a unique vector $v^k_{ij}:=v_0(\tilde{\Gamma}^k_{ij})^{\hconj} \in\C^N$ such that $\Gamma^k_{ij}=\tilde{\Gamma}^k_{ij}p+\delta^k_j\der_i(p)=(v^k_{ij})^{\hconj}v_0+0$, leading to the characterization of $\nabla$ being simplified through the use of the standard scalar product on $\C^N$:
\begin{equation*}
    e_k\Gamma^k_{ij}=(0,...,\alpha_kv_0,0,...,0)\cdot (v^k_{ij})^{\hconj}v_0=(0,...,\alpha_k (v_0 (v^k_{ij})^{\hconj})v_0,0,...,0)=e_k(v_0 (v^k_{ij})^{\hconj}),
\end{equation*}
and if we let $\lambda^k_{ij}=v_0 (v^k_{ij})^{\hconj}\in\C$ we have that
\begin{equation*}
    \nabla_i e_j=e_k\Gamma^k_{ij}=e_k\lambda^k_{ij};
\end{equation*}
since the matrix $\tilde{\Gamma}^k_{ij}$ that yields $\lambda^k_{ij}$ is arbitrary, the above expression defines an affine connection for any choice of $\lambda^k_{ij}\in\C$.

With this simplification, the hermiticity condition for $h(e_i,\nabla_j e_k)$ becomes
\begin{equation*}
    h(e_i,\nabla_j e_k)=(\tilde{\mu}_{il}\lambda^l_{jk})v_0^{\hconj}v_0=(\tilde{\mu}_{il}\bar{\lambda}^l_{jk})v_0^{\hconj}v_0=h(e_i,\nabla_j e_k)^{\hconj},
\end{equation*}
and multiplying by $\tilde{\mu}^{mi}$ from the left, we see that this is true if and only if $\lambda^l_{jk}\in\R$ for every $j,k,l$.

Since the torsion condition is not trivial if dim $\g>1$ we may use Koszul's formula to find constants $\lambda^k_{ij}$ satisfying the metric and torsion-free conditions for $\nabla$ directly if they exist. We find that the connection must satisfy
\begin{multline*}
    2h(\nabla_i e_j,e_k)=\der_i h(e_j,e_k)+\der_j h(e_i,e_k)-\der_k h(e_i,e_j)\\
    -h(e_i,\varphi([\der_j,\der_k]))+h(e_j,\varphi([\der_k,\der_i]))+h(e_k,\varphi([\der_i,\der_j]))=0,
\end{multline*}
since the Lie bracket is zero and $\der_i (h_{jk})=0$ for every $i,j,k$. 
On the other hand, $h(\nabla_i e_j,e_k)$ may be calculated explicitly explicitly:
\begin{equation*}
    0=2h(\nabla_i e_j,e_k)=2h(e_l\lambda^l_{ij},e_k)=2\tilde{\mu}_{lk}\bar{\lambda}^l_{ij};
\end{equation*}
multiplying by $\tilde{\mu}^{mk}$ from the right, one finds that this happens if and only if $\bar{\lambda}^k_{ij}=\lambda^k_{ij}=0$ for every $i,j,k$. In other words, for an arbitrary vector $v=e_kB^k\in (\C^N)^n$ and $\der\in\g$,
\begin{equation*}
    \nabla_{\der} v=\nabla_{\der} (e_k B^k)=e_k\der(B^k)=e_k[\hat{D},B^k]=\lambda_{\der} e_k B^k-(e_k B^k)\cdot \hat{D}=\lambda_{\der} v-v\cdot \hat{D},
\end{equation*}
where $\hat{D}=\hat{D}(\der)$ and $\lambda_{\der}$ is the eigenvalue of $\hat{D}$ associated with $v_0$. We summarize the above discussion as follows.

\begin{proposition}
    Let $\g$ be an abelian Lie algebra with basis $\{\der_1,...,\der_n\}$ and  let $(C_{\A},h)$ be a real metric calculus such that
    \begin{equation*}
        C_{\A}=(\Mat{N},\g_D,(\C^N)^n,\varphi)
    \end{equation*}
    and $\varphi(\der_i)=(0,...,0,\alpha_i v_0,0,...,0)\in(\C^{N})^n$,
    where $\alpha_i\in\Rstar$ and $v_0$ is a unit vector.
    Then there exists a unique connection $\nabla$ such that $(C_{\A},h,\nabla)$ is pseudo-Riemannian if and only if $v_0$ is an eigenvector of $\hat{D}(\der_i)$ with eigenvalue $\lambda_i$ for $i=1,...,n$. In this case, the Levi-Civita connection $\nabla$ is given by
    \begin{equation*}
        \nabla_{\der_i} v=\lambda_i v-v\cdot\hat{D}(\der_i),
    \end{equation*}
    for $v\in (\C^N)^n$.
\end{proposition}

\bibliographystyle{alpha} 
\bibliography{references} 

\end{document}